\documentclass[12pt]{article}  
\usepackage{amssymb}
\usepackage{amsmath,amsthm,enumerate}
\theoremstyle{plain}
\newtheorem{thm}{Theorem}[section]
\newtheorem{prop}[thm]{Proposition}
\newtheorem{cor}[thm]{Corollary}
\newtheorem{lemma}[thm]{Lemma}

\theoremstyle{definition}

\newtheorem{rmk}[thm]{Remark}

\newenvironment{enum1}{\begin{enumerate}[\upshape (1)]}{\end{enumerate}}

\title{Sufficient and equivalent criteria for the Riemann Hypothesis\thanks{Research
    supported by Swiss National Science Foundation Grant no. 107887. Support has also been given over time by Scuola Normale Superiore of Pisa, University of Z\"urich and IBM Z\"urich Research Lab.}}
\author{Davide Schipani\\
Institute of Mathematics\\
University of Zurich\\
davide.schipani(at)math.uzh.ch}
\date{\today}
\begin{document}\maketitle

\begin{abstract}


The paper presents several new sufficient conditions, as well as new equivalent criteria for the classical Riemann Hypothesis. 
Noteworthy are also other statements and remarks about $\zeta$ to be found throughout the paper.

\end{abstract}\medskip

\section{Introduction}

It is well known that the Riemann zeta function $\zeta(s)$ is zero at $s=-2,-4,-6,\ldots$, which are called the trivial zeros. The Riemann Hypothesis conjectures that all non-trivial zeros have real part equal to $1/2$. This conjecture inspired many mathematicians, so that a vast literature is nowadays available, and it has been also computationally verified for the first $10^{13}$ non trivial zeros, ordered by increasing positive imaginary part (see e.g. \cite{go04}). Several authors also came up with sufficient and/or necessary conditions: for many of them, and at the same time to have a more general overview of the problem, we can refer the reader for example to \cite{bo06}, \cite{bo08}, \cite{co03} and the references therein. This work is mainly devoted to present new sufficient conditions and equivalent criteria.

An overview of the paper is the following. After recalling below some properties of the zeta function, 
we will state and prove in section $2$ sufficient conditions for the Riemann Hypothesis while in section $3$ we deal with equivalent criteria. 

In a follow-up paper we will state analogues of some of these results for generalized versions of the Riemannn Hypothesis (\cite{sc10}).


Throughout this document we make use of basic properties of the Riemann zeta function: we summarize them below so that also a non-specialist can follow the arguments. We will refer later to them by using the numbers in round brackets. Some facts about the $\Gamma$ function are also added.

\begin{enum1}
\item It is known that there are no nontrivial zeros of $\zeta$ in the half-planes $x\leq 0$ and $x\geq 1$ (nor in some regions in between, for example on the real axis between $0$ and $1$: see e.g. \cite[chapter 13]{ap76} or \cite{sc10b}). So
  the Riemann Hypothesis says essentially that the only zeros in the
  critical strip, i.e in the strip $0<\Re(s)<1$, are on the line
  $x=\frac{1}{2}$. The critical strip $0<\Re(s)<1$ will then be the
  region of interest for us, in particular as far as
  analytical properties or
  convergence of sequences of functions are concerned.
\item The function
  $\Lambda(s)=\pi^{-\frac{s}{2}}\Gamma(\frac{s}{2})\zeta(s)$ is
  analytic in the whole complex plane apart from poles at $0$ and
  $1$. Furthermore $\Lambda(s)=\Lambda(1-s)$ (see e.g \cite{la06}).
\item Zeros of $\zeta$ in the critical strip are the same as
  zeros of $\Lambda$ ($\Gamma$ has no zeros). 
  It follows from $(2)$ that if $s$ is a zero of
  $\zeta$, then also $1-s$ is a zero.
\item Since $\zeta$ is analytic in the critical strip, its zeros
  are isolated.
\item A representation for $\zeta$ in the strip is the following:
$$
\zeta(s)=(1-2^{1-s})^{-1}\cdot  \phi(s),
$$
 where $\phi(s)=\sum_{n=1}^\infty  \frac{(-1)^{n-1}}{n^s}$ (see e.g. \cite[section 10.2]{na00}). Again, zeros of $\zeta$ in the strip are the same as those of $\phi$.
\item $\overline{\zeta(s)}=\zeta(\bar{s})$. (It follows easily, for example,
  by the representation in $(5)$
; or one can see it from
  the Schwarz reflection principle in complex analysis (for example \cite{ca95}), since
  $\zeta$ is real on the real axis.)
\item From $(3)$ and $(6)$ it follows that zeros are symmetric both with
  respect to the critical line $x=\frac{1}{2}$ and with respect to the real axis.
\item By Euler-Maclaurin summation formula (see e.g \cite[section 6.4]{ed01}), if $s=\sigma+it$ belongs to $\mathbb{C}\backslash\{1\}$ the following holds:
$$
\zeta(s)=\sum_{n=1}^{N-1}n^{-s}+\frac{N^{1-s}}{s-1}+\frac{1}{2}N^{-s}+O(N^{-\sigma-1});
$$
Moreover, by the upper bound on the remainder terms as expressed in
\cite{ed01}, the constant implicit in the $O$-notation can be chosen to be the
same for all $s$ in any compact subset of $\mathbb{C}\backslash\{1\}$.

Actually, if we consider a region $\sigma\geq \sigma_0>0$, with $\sigma_0$ fixed, we could even say that
$$
\zeta(s)=\sum_{n=1}^{N}n^{-s}+\frac{N^{1-s}}{s-1}+O(N^{-\sigma})
$$
holds uniformly, provided that $N>C|t|/2\pi$, where $C$ is a given constant greater than $1$ (see \cite[Theorem 4.11]{ti86}).

\item $\overline{\Gamma(s)}=\Gamma(\bar{s})$. (It follows for example from
  the definition of $\Gamma$ as a limit 
 or again by the Schwarz reflection principle.)
\item If $-1/2\leq\sigma\leq 1/2$, then 
  $\left|\frac{\Gamma(\frac{1-s}{2})}{\Gamma(\frac{s}{2})}\right|\leq|(1+s)/2|^{1/2-\sigma}$ (\cite[Lemma 1] {ra59} specialized with $q=0$).


\end{enum1}

\section{Sufficient conditions for the Riemann Hypothesis}

Before stating our first sufficient condition, we recall a result of Selberg (\cite{se42}, \cite[section 2.1.11]{ka95}) which says that 

$$
\lim_{T\to\infty} \frac{1}{T}\mu\left\{\tau\in(0,T] : \forall k\in(\tau,\tau+\Phi(\tau)/\log \tau)
\ \zeta(\frac{1}{2}+ik)\neq 0   \right\}=0,
$$
where $\Phi(\tau)$ is any positive function which tends to $\infty$ as $\tau\to\infty$ and $\mu$ stands for the Lebesgue measure.

By taking for example $\Phi(\tau)=\log\log \tau$, this implies that for any fixed $\eta$
$$
\lim_{T\to\infty} \frac{1}{T}\mu\left\{\tau\in(0,T] : \forall k\in(\tau,\tau+\eta)
\ \zeta(\frac{1}{2}+ik)\neq 0   \right\}=0,
$$ 
or
$$
\lim_{T\to\infty} \frac{1}{T}\mu\left\{\tau\in(0,T] : \exists k\in(\tau,\tau+\eta) \mid
\zeta(\frac{1}{2}+ik)=0 \right\}=1.
$$

We will also need another important result, namely Voronin's universality theorem, which says the following (see e.g. \cite{st07}): suppose that $K$ is a compact subset of the strip
$\frac{1}{2}<\sigma<1$ with connected complement and let $g(s)$
be a non-vanishing continuous function on $K$ which is analytic
in the interior of $K$; then 
for any
$\epsilon>0$
$$
\liminf_{T\to\infty} \frac{1}{T} \mu\left\{\tau\in[0,T] :
\max_{s\in K} |\zeta(s+i\tau)-g(s)|<\epsilon   \right\}>0.
$$

We remark (see \cite[section 8.1]{st07}) that, if Voronin's theorem were true
even if $g(s)$ is allowed to vanish, then the Riemann Hypothesis would be
false. Though, as shown in the reference cited, this cannot
happen, i.e a function having zeros cannot be approximated
uniformly by $\zeta$, 
 which
in fact 
hints at the relation between universality and distribution
of zeros. In this regard we notice that with another notion of universality, slightly different from that expressed by Voronin's theorem, it would be possible to find functions
which satisfy the main properties of $\zeta$, namely the
functional equation in $(2)$, the property of being analytic except for a
pole in $1$ and that of being real on the real axis, but do not
satisfy the Riemann Hypothesis (see \cite{ni09}).

Furthermore, not only the zeta function is universal (in the sense of the formula above), but its derivative is even strongly universal in Voronin's sense, strongly meaning that the function to be approximated needs not to be non-vanishing; see \cite[section 1.3 and 1.6]{st07}.

We are now ready for our first statement:

\begin{thm}\label{suff2}
Suppose that $\zeta'$ is strongly universal for any compact subset of the region $\{\frac{1}{2}<\sigma<1\}\cup\{[s_*,s_*+i\eta]\}$, where $s_*$ is a zero on the critical line and $\eta$ is a positive real number. Then the Riemann Hypothesis is true.
\end{thm}


\begin{proof}

We are going to proceed by contradiction: we assume that the Riemann Hypothesis is false, which implies by $(7)$ that there is a zero, $s_0$, in the strip
$\frac{1}{2}<\sigma<1$, and let $K$ be a compact subset in this
strip containing $s_0$. 

 For any given
$\epsilon>0$, by continuity, there exists $\eta>0$ such that
$|s_*-s|<\eta$ implies $|\zeta(s)|<\epsilon/2$. 


Because of the hypothesis, for any $\delta>0$
\[
\liminf_{T\to\infty} \frac{1}{T} \mu\left\{\tau\in[0,T] :
\max_{s\in K_2} |\zeta'(s+i\tau)-\zeta'(s)|\leq \delta   \right\}>0,
\]
 
for any compact $K_2$ in $\{\frac{1}{2}<\sigma<1\}\cup\{[s_*,s_*+i\eta]\}$.


We will take $K_2$ as a compact
containing all lines from $K$ to the segment $[s_*,s_*+i\eta]$.
 



Let also $\delta$ be equal to $\frac{\epsilon}{2(\max_{s\in
    K}|s-s_*|+\eta)}$ and consider $|\zeta(s+i\tau)-\zeta(s)|$ for each of those $\tau$ such that $\max_{s\in K_2} |\zeta'(s+i\tau)-\zeta'(s)|\leq \delta$.

By the triangular inequality we have for any $0\leq\alpha\leq 1$ and $s\in K$
$$
|\zeta(s+i\tau)-\zeta(s)|\leq\left|\int_{s_*+i\alpha\eta}^s(\zeta'(s+i\tau)-\zeta'(s))ds\right|+|\zeta(s_*+i\alpha\eta+i\tau)|+|\zeta(s_*+i\alpha\eta)|
$$
where the path of integration is chosen to be a straight line connecting the two points, so that the first part on the right can be bounded by $|s-(s_*+i\alpha\eta)|\cdot\delta$, which is less or equal than $\epsilon/2$.

Now 
 we know by Selberg's theorem that for almost all $\tau$ there exists $0\leq\alpha\leq
1$ such that $\zeta(s_*+i\alpha\eta+i\tau)=0$; moreover $|\zeta(s_*+i\alpha\eta)|\leq \epsilon/2$ by the choice of $\eta$, so that the second part is less or equal than $\epsilon/2$, too.

Thus in general we can conclude that
$$
\liminf_{T\to\infty} \frac{1}{T} \mu\left\{\tau\in[0,T] :
\max_{s\in K} |\zeta(s+i\tau)-\zeta(s)|\leq \epsilon   \right\}>0.
$$

Though, as remarked before stating the theorem, this cannot happen, since $\zeta$ is
supposed to vanish on $K$, and we obtain a contradiction.

\end{proof}

\begin{rmk}
Alternatively, instead of proving it by contradiction, we could
have used directly a theorem by Bagchi (see e.g. \cite[section 8.2]{st07}), which says that the Riemann Hypothesis is true if and only if for any compact subset $K$ of the strip
$\frac{1}{2}<\sigma<1$ with connected complement and for any
$\epsilon>0$
$$
\liminf_{T\to\infty} \frac{1}{T} \mu\left\{\tau\in[0,T] :
\max_{s\in K} |\zeta(s+i\tau)-\zeta(s)|<\epsilon   \right\}>0.
$$

\end{rmk}

\begin{rmk}
By looking at the given proof one may think about varying it by considering the possibility that zeta takes a particular value on a line $\sigma=\sigma_0$, with $1/2<\sigma_0<1$, as often as it happens for zeros on the critical line. Though it has been shown that for any complex number $c\neq 0$, the number of roots of $\zeta(s)=c$ up to height $T$ in the strip $1/2<\sigma_1<\sigma_2<1$ is only asymptotic to $KT$ for a finite positive constant $K$ (\cite{bo32}, \cite[Theorem 1.5]{st07}).
\end{rmk}

\begin{cor}\label{coro}
Suppose that $\zeta'$ is strongly universal in $\frac{1}{2}\leq\sigma<1$. Then the Riemann Hypothesis is true.
\end{cor}



\vspace{1cm}

We will find now another sufficient condition, by a completely different approach.

We first define:

$$
F(s)\doteq \frac{\zeta(s)}{\zeta(1-s)}.
$$ 

Notice that $F(s)$ can be
analytically continued in the zeros of $\zeta$, by defining it
there to
be $F(s)=\frac{\pi^{s-\frac{1}{2}}\Gamma(\frac{1-s}{2})}{\Gamma(\frac{s}{2})}$
from $(2)$. 

This actually not only shows that zeros of the denominator
are just removable singularities, but also that if $s$ is a
zero, then $s$ and $1-s$ have the same multiplicity. \footnote{Just suppose that it is not the
case, then $s$ has higher multiplicity, in order for $F(s)$ to
remain bounded in a neighboorhood of $s$ (and thus to be
extendable: see \cite[section III.4.4]{ca95}); but in this case
$\frac{\zeta(1-s)}{\zeta(s)}$ would not be extendable, by the
same reason. Or just notice that if $s$ is a zero with a higher multiplicity
than $1-s$, then $F(s)=0$, while from the explicit expression above we see that $F(s)$ is never zero
in the critical strip.}

Consider now $H_N(s)\doteq\sum_{n=1}^{N}n^{-s}+\frac{N^{1-s}}{s-1}$, which can also be seen as a truncated sum: if we define (see also \cite{pu08}):
$$
h_n(s)=
\begin{cases}
\frac{1}{n^s}-\frac{n^{1-s}-(n-1)^{1-s}}{1-s}\ \ \text{if}\ \  n\geq 2\\
1+\frac{1}{s-1}\ \  \text{if}\ \  n=1,
\end{cases}
$$
then $H_N(s)$ is just $\sum_{n=1}^N
  h_n(s)$, the $N$-th partial sum of
  the series $\sum_{n=1}^{\infty} h_n(s)$. Notice that because of the summation formula in $(8)$ we have:

$$
\zeta(s)=\sum_{n=1}^{N} h_n(s)-\frac{1}{2}N^{-s}+O(N^{-\sigma-1})
$$
so that 
$$
\zeta(s)=\sum_{n=1}^{\infty} h_n(s)
$$
and if $\zeta(s)=0$,
$$
\sum_{n=1}^{N} h_n(s)=\frac{1}{2}N^{-s}+O(N^{-\sigma-1}).
$$


Suppose now $s_*=\sigma_* +it_*$ is in the
  critical strip.


\begin{prop}\label{prop1}
Whenever $\zeta(s_*)\neq 0$ (and so $\zeta(1-s_*)\neq 0$), then
$$\lim_{s \to s_*} \frac{\zeta(s)}{\zeta(1-s)}=\lim_{N \to
  \infty} \frac{H_N(s_*)}{H_N(1-s_*)}.
$$
\end{prop}
\begin{proof}
Since $\lim_{N \to
  \infty} H_N(1-s_*)=\zeta(1-s_*)\neq 0$, 
we have: 
$$
\lim_{N \to
  \infty} \frac{H_N(s_*)}{H_N(1-s_*)}=\frac{\lim_{N \to
  \infty} H_N(s_*)}{\lim_{N \to
  \infty} H_N(1-s_*)}=
\frac{\zeta(s_*)}{\zeta(1-s_*)}=F(s_*).
$$
\end{proof}
Now, if $a(s)$ and $b(s)$ are functions with real positive values, we define $f_{a(s)}(N,s)$ as $\lceil a(s)^{-1/(2\sigma)}
N^{1-\sigma} \rceil$ (resp. with $b$), 
with, as usual, $s=\sigma +it$, and let
$$
H^{a,b}(s)\doteq\lim_{N \to \infty} \left|\frac{H_{f_{a(s)}(N,s)}(s)}{H_{f_{b(1-s)}(N,1-s)}(1-s)}\right|.
$$


Notice that in proving the previous proposition we didn't use the fact
that we had $N$-th partial sums in both the numerator and the
denominator. Then along the same lines we can also prove:

\begin{prop}\label{prop2}
If $\zeta(s_*)\neq 0$, then
$$|F(s_*)|=\lim_{s \to s_*} \left|\frac{\zeta(s)}{\zeta(1-s)}\right|=\lim_{N \to \infty} \left|\frac{H_{f_{a(s_*)}(N,s_*)}(s_*)}{H_{f_{b(1-s_*)}(N,1-s_*)}(1-s_*)}\right|=H^{a,b}(s_*).
$$
\end{prop}

Consider now $H^{a,b}(s_*)$ for $s_*$ being a zero of $\zeta$: we
claim that 
\begin{lemma}\label{lemma1}
For any zero $s_*$ of $\zeta$, if $a(s),b(s)$ are constants $A$ and $B$, we have that $H^{a,b}(s_*)=|\sqrt{A}/\sqrt{B}|$. 
\end{lemma}
\begin{proof}
Write $\lceil A^{-1/(2\sigma)}
N^{1-\sigma} \rceil$  as $A^{-1/(2\sigma)}
N^{1-\sigma}+\epsilon(N)$ with $0\leq\epsilon(N)<1$, and then as $(A^{-1/(2\sigma)}
N^{1-\sigma})\cdot (1+\frac{\epsilon(N)}{A^{-1/(2\sigma)}
N^{1-\sigma}})$. We do a similar decomposition for
$f_B(N,1-s)$, 
and since we are supposing that $\zeta(s_*)=0$ (so that
$\sum_{n=1}^{N} h_n(s_*)=\frac{1}{2}N^{-s_*}+O(N^{-\sigma_*-1})$
as shown above using the summation formula),
we can write
$$
H^{a,b}(s_*)=\lim_{N\to\infty}\left|\frac{\frac{1}{2}(A^{-1/(2\sigma_*)}
N^{1-\sigma_*})^{-s_*}(1+\frac{\epsilon_1(N)}{A^{-1/(2\sigma_*)}
N^{1-\sigma_*}})^{-s_*}+O(f_A(N,s_*)^{-\sigma_*-1})}{\frac{1}{2}(B^{-1/(2-2\sigma_*)}N^{\sigma_*})^{-1+s_*}(1+\frac{\epsilon_2(N)}{
B^{-1/(2-2\sigma_*)}N^{\sigma_*}})^{-1+s_*}+O(f_B(N,1-s_*)^{-2+\sigma_*})}\right|.
$$
By the binomial expansion of $(1+x)^z$,
$$
(1+\frac{\epsilon_1(N)}{A^{-1/(2\sigma_*)}
  N^{1-\sigma_*}})^{-s_*}=1+O(\frac{1}{N^{1-\sigma_*}})
$$ 
and
$$
(1+\frac{\epsilon_2(N)}{B^{-1/(2-2\sigma_*)}N^{\sigma_*}})^{-1+s_*}=1+O(\frac{1}{N^{\sigma_*}}),
$$
so that we
finally get
\[
\lim_{N\to\infty}\left|\frac{\frac{1}{2}A^{s_*/(2\sigma_*)}
N^{-s_*(1-\sigma_*)}+O((N^{1-\sigma_*})^{-\sigma_*-1})}{\frac{1}{2}B^{(1-s_*)/(2-2\sigma_*)}
N^{(-1+s_*)\sigma_*}+O((N^{\sigma_*})^{-2+\sigma_*})}\right|=
\]
\[
=\lim_{N\to\infty}\left|\frac{\frac{1}{2}A^{s_*/(2\sigma_*)}
N^{-s_*(1-\sigma_*)}}{\frac{1}{2}B^{(1-s_*)/(2-2\sigma_*)}
N^{(-1+s_*)\sigma_*}}\right|=\lim_{N\to\infty}\left|\frac{A^{\sigma_*/(2\sigma_*)}
N^{-\sigma_*(1-\sigma_*)}}{B^{(1-\sigma_*)/(2-2\sigma_*)}
N^{(-1+\sigma_*)\sigma_*}}\right|=|\sqrt{A}/\sqrt{B}|
\]

\end{proof}

Consider now the case $a(s)=b(s)=|F(s)|$. Using the
fact that $|F(s)|=1/|F(1-s)|$, one sees in a similar manner that
$$
H^{|F(s)|,|F(s)|}(s)=|\sqrt{|F(s)|}/\sqrt{|F(1-s)|}|=|F(s)|
$$
if $s$ is a zero. If it is not a zero, then the statement is of course true, so that we have built a continuous function.

We look at this point at the functions which are the argument of $\lim$ in the definition of $H^{|F(s)|,|F(s)|}(s)$, and call them $|F(s)|_N$. We first notice that for $\sigma=1/2$, they are always equal to $1$ and equal to the limit $|F(s)|=1$. In the right open half instead we have $|F(s)|<1$ as soon as $t$ is large enough ($\left|\frac{\zeta(\sigma+it)}{\zeta(1-\sigma-it)}\right|$ is actually monotonically decreasing in $\sigma$ for fixed $t>2\pi+1$, as proven in \cite[Theorem 1]{sa03}). 
Moreover it is not hard to see that, for each $s$ in this region, the sequence $\{|F(s)|_N\}_N$ admits a strictly monotone subsequence: for that it suffices to check that $\{|F(s)|_N\}_N$ does not get constant from some point, which can be seen for example by the fact that, for any two constants $k_1$ and $k_2$, $\lceil k_1
N^{1-\sigma} \rceil$ and $\lceil k_2
N^{\sigma} \rceil$ are increasing by different rates for $N$ large enough. Now, if for an $s$ in the half strip $|F(s)|_N<|F(s)|_M$ for some $N,M$, by continuity this holds true also in a neighbourhood. In the theorem below we will essentially assume that the intersection of all these kinds of neighbourhoods relative to some subsequence does not reduce to a point.

Before stating the theorem, we notice that we can also consider slightly different functions converging to $|F(s)|$: for example in the numerator of $|F(s)|_N$, we could take $m_{|F(s)|}(N,s)\doteq 1+\lceil |F(s)|^{-1/2\sigma}
N^{1-\sigma} \rceil$ instead of $f_{|F(s)|}(N,s)$ and call the new approximants $|F_m(s)|_N$.

\begin{thm}
Suppose that for each $s$ in the right (or left) open half strip we can find a neighbourhood $I_s$ and a subsequence $N_k(s)$ such that $\{|F(s)|_{N_k}\}_{N_k}$ and $\{|F_m(s)|_{N_k}\}_{N_k}$ converge monotonically (either increasing or decreasing for the whole $I_s$). Then the Riemann Hypothesis is true.
\end{thm}
\begin{proof}

Thanks to the hypothesis of monotonicity, and since the approximants converge by construction to a continuous function, by Dini's theorem (see e.g. \cite[Theorem 7.13]{ru76}),
 the convergence in a compact set $D$ inside $I_s$ would be uniform. This implies that, for $N$ bigger than some $N_1$,  $|F(s)|_N$ is uniformly bounded in $D$: then, for those $N$, possible zeros of the denominator $H_{f_{|F(1-s)|}(N,1-s)}(1-s)$ would have to be zeros of the numerator $H_{f_{|F(s)|}(N,s)}(s)$, too. 
And for $N$ bigger than some $N_2$ zeros of $H_{f_{|F(1-s)|}(N,1-s)}(1-s)$ would also be zeros of $H_{m_{|F(s)|}(N,s)}(s)$. But $H_N$ and $H_{N+1}$ cannot both be zero, at least for $N$ large, which is what interests us now: the difference between the two is $h_N(s)$, which is $0$ exactly when
$1-s=n-(n-1)^{1-s}n^s=n-n(1-\frac{1}{n})^{1-s}$; though
$1-s=n-n(1-\frac{1-s}{n})$ which is only an approximation of the
binomial expansion of $n-n(1-\frac{1}{n})^{1-s}$. (In fact in the critical strip this is also true for small values, as stated in \cite{pu04}.)

In conclusion, for all $N$ big enough, $H_N$ would have to be nonzero throughout $D$, and this is sufficient to prove that $\zeta$ has no zeros in $D$ by a theorem of Hurwitz in complex analysis. Hurwitz proved in fact the following, as stated in \cite[Theorem 5.13]{st07}: 
let $G$ be a region and $\{f_n\}$ be a sequence of functions
analytic on $G$ which converges uniformly on $G$ to some function
$f$; suppose also that $f(s)$ is not identically zero, then an
interior point $s_0$ of $G$ is a zero of $f(s)$ if and only if
there exists a sequence $\{s_n\}$ in $G$ which tends to $s_0$ as
$n\to \infty$, and $f_n(s_n)=0$ for all sufficiently large $n$.

\end{proof}



\section{Equivalent criteria for the Riemann Hypothesis}

In the first part of this section we concentrate on results that involve zeros of partial sums. One of the most famous theorems in this regard  is a classical
by Tur\'an, who proved the following (see e.g. \cite[section 8.14]{ap90} or \cite[section 8.3]{bo08}): let $\zeta_N(s)\doteq\sum_{n=1}^N \frac{1}{n^s}$; if there exists
an $N_0$ such that $\zeta_N(s)\neq 0$ for all $N\geq N_0$ and all
$\sigma>1$, then the Riemann Hypothesis is true.

We remark that this is only a one-side implication, which has also
proved to be unuseful for a proof of the Riemann Hypothesis, since Montgomery proved
the existence of roots of these partial sums in certain half planes to
the right of the critical strip, as cited in the references
above. 

In this regard and in general about zeros of partial sums it's
recommended to have a look at \cite{bo07}, \cite{go07b}, \cite{go08a} too.

Our first result in this direction is the following:
\begin{thm}
The Riemann Hypothesis is true if and only if for any compact disc $K$ in the right (or left)
open half of the critical strip there exists an $N_0$ such
that for infinitely many $N>N_0$, $H_N(s)$ is non-vanishing for all $s$ in
$K$.
\end{thm}
\begin{proof}
If the Riemann Hypothesis is true, then
there are no zeros in the right open half. In any compact $K$
then $|\zeta|$ has a minimum $m>0$. Now, by $(8)$, we know that
$$
|H_N(s)-\zeta(s)|=|\frac{1}{2}N^{-s}+O(N^{-\sigma-1})|
$$
(keeping in mind the remark about the $O$ notation at $(8)$).
By taking $N$ big enough, we can make it smaller than $m/2$ for all $s$ in $K$, so that, by the triangular inequality,
$|H_N(s)|=|H_N(s)-\zeta(s)+\zeta(s)|\geq |\zeta(s)|-|H_N(s)-\zeta(s)|\geq m/2>0$.

In the other direction, to prove that the Riemann Hypothesis is
true, we can use Hurwitz's theorem (see previous section). Since our
hypothesis implies that for any compact disc $K$
there cannot exist a sequence $s_N$ in $K$ which is a zero for
$H_N$ for all sufficiently large $N$, Hurwitz's theorem tells us
that there is no zero of $\zeta$ in $K$. But any $s$ in the right
open half is center of a compact disc in the interior of this strip.
 
\end{proof}

We can derive a similar statement using the approximants $\phi_N(s)\doteq\sum_{n=1}^N  \frac{(-1)^{n-1}}{n^s}$:
\begin{thm}
The Riemann Hypothesis is true if and only if for any compact disc $K$ in the right
open half of the critical strip there exists an $N_0$ such
that, for infinitely many $N>N_0$, $\phi_N(s)$ is non-vanishing for all $s$ in
$K$.
\end{thm}
\begin{proof}
We can follow the same argument as above, provided that we consider $(1-2^{1-s})^{-1}\cdot \phi_N(s)$. Notice that in this case we know that $\phi_N(s)$ is converging uniformly in any compact of the strip, since any Dirichlet series converges uniformly on every compact subset interior to its half-plane of convergence (see e.g. \cite[chapter 11]{ap76}).
\end{proof}

Generalizing a bit these results, for example for the case of $H_N(s)$:
\begin{thm}
The Riemann Hypothesis is true if and only if for any compact disc $K$ in the right
open half of the critical strip there exists an analytic function
$f(N,s)$
which tends to zero uniformly in $K$ as $N\to\infty$, $\alpha_1,\ldots,\alpha_{d}\in\mathbb{C}$, $L_1,\ldots,L_d,N_0\in\mathbb{N}$ such
that, for infinitely many $N>N_0$, $f(N,s)+\alpha_1 H_N(s)+\alpha_2
H_{N+L_1}+\cdots+\alpha_d
H_{N+L_{d-1}}+(1-\sum_{i=1}^d \alpha_i) H_{N+L_d}$ is non-vanishing for all $s$ in
$K$.
\end{thm}
\begin{proof}
If the Riemann Hypothesis is true, then for any $K$, as we saw in the previous equivalence, there
exists an $N_0$ such
that, for infinitely many $N>N_0$, $H_N(s)$ is non-vanishing for all $s$ in
$K$. This corresponds to taking $f(N,s)=0$ for all $s$ in $K$, $\alpha_1=1$ and all
other $\alpha$ equal to $0$.
In the other direction, the proof follows the one given above for
the less general case, since
$f(N,s)+\alpha_1 H_N(s)+\alpha_2
H_{N+L_1}+\cdots+(1-\sum_{i=1}^d \alpha_i) H_{N+L_d}$ is also
converging to $\zeta(s)$ uniformly on
$K$.
\end{proof}

For
example it would be sufficient to prove that the arithmetic mean
of $H_N$ and $H_{N-1}$ is always nonzero for infinitely many $N>N_0$ and for
all $s$ in any given $K$. This might be easier, though being
essentially equivalent; we already remarked that if one of the two is zero, the other is not, which might be useful.

\vspace{1cm}

In the sequel other equivalent reformulations, though of a different nature, are to be found.

We first notice that, if $s_0=\sigma_0+it_0$ is a zero of $\zeta$
in the critical strip, then
$F(s_0)$, where $F$ is defined as in the previous section, is a nonzero value, namely
$F(s_0)=\frac{\pi^{s_0-\frac{1}{2}}\Gamma(\frac{1-s_0}{2})}{\Gamma(\frac{s_0}{2})}$.

\begin{thm}\label{teo1}
The Riemann Hypothesis is true if and only if $|F(s_0)|=1$ for any zero $s_0$ of $\zeta$.
\end{thm}
\begin{proof}

 
When $\sigma=1/2$, $|F(s)|$ is equal to $1$ by $(6)$.

On the other hand, if $\sigma\neq 1/2$ and $t\geq 2\pi+1$, by \cite[Theorem 1]{sa03}, we have $|F(s)|\neq 1$. (Notice that thanks to $(6)$ and the fact that $F(s)=1/F(1-s)$, if $|F(s)|\neq 1$ in a point, neither it will be in its symmetric with respect to $1/2$ or in its symmetric with respect to the critical line.)

But the region $\{t<2\pi+1\}$ has no influence in this context, since we know that the first zero in the upper half-plane has imaginary part larger than $14$.

\end{proof}

\begin{rmk}
We can actually almost complement the result in \cite{sa03} mentioned in the proof, improving somehow the known asympotic estimates or almost equivalences (see for example, besides \cite{sa03}, also \cite[Lemma 6.1]{go07b} and \cite[Theorem 4.3]{br04} in addition to the known asymptotic estimate by Stirling's approximation). 

First notice that the proof in \cite{sa03} is derived ultimately through some estimates and actually requires $t>t_0$ where $t_0$ is such that $\left|\frac{1}{2}+it_0\right|=2\pi e^{0.0212411}=6.418\ldots$.

We take care instead of the region below and prove the following:
\begin{prop}\label{teo2}
$|F(s)|\neq 1$ in the regions $\{0<\sigma<1/2\}\cap\{\sqrt{(1+\sigma)^2+t^2}<2\pi\}$ or $\{1/2<\sigma<1\}\cap\{\sqrt{(2-\sigma)^2+t^2}<
2\pi\}$. 
\end{prop}
\begin{proof}
In fact $|F(s)|=1$
happens when
$\pi^{\sigma-1/2}=|\pi^{s-1/2}|=\left|\frac{\Gamma(\frac{s}{2})}{\Gamma(\frac{1-s}{2})}\right|$.

From ($10$) 
we get that, for $\sigma<1/2$,
$\left|\frac{\Gamma(\frac{s}{2})}{\Gamma(\frac{1-s}{2})}\right|\geq
|(1+s)/2|^{\sigma-1/2}$. This means that, for $\sigma<1/2$ and
$\sqrt{(1+\sigma)^2+t^2}< 2\pi$, $\left|\frac{\Gamma(\frac{s}{2})}{\Gamma(\frac{1-s}{2})}\right|>
\pi^{\sigma-1/2}$, so that the right
and left member will never be equal in this case.  

For $\sigma>1/2$, as pointed out before, we just need to consider the symmetric region with respect to the critical line. Or, equivalently, we let $z=1-s$ so that $\Re{z}<1/2$: then $\left|\frac{\Gamma(\frac{s}{2})}{\Gamma(\frac{1-s}{2})}\right|=\left|\frac{\Gamma(\frac{1-z}{2})}{\Gamma(\frac{z}{2})}\right|\leq|(1+z)/2|^{1/2-\Re(z)}=|(2-s)/2|^{\sigma-1/2}$. And this
means that, for $\sigma>1/2$ and $\sqrt{(2-\sigma)^2+t^2}< 2\pi$, we have that $\left|\frac{\Gamma(\frac{s}{2})}{\Gamma(\frac{1-s}{2})}\right|<
\pi^{\sigma-1/2}$, and again the right
and left member will never be equal in this case either.
\end{proof}

It actually follows from the proof above that, in $\{0<\sigma<1/2\}\cap\{\sqrt{(1+\sigma)^2+t^2}<2\pi\}$,  $|F(s)|$ is less than $1$ (or, alternatively, since $|F(s)|$ is continuous, it is enough for example to consider the part of real axis inside the strip, where $\zeta$ is negative and decreasing: see \cite[notes to chapter 1]{iv03}). Instead, in $\{0<\sigma<1/2\}\cap\{t>t_0\}$, $|F(s)|$ is bigger than $1$ by \cite[Theorem 1]{sa03}.

We remark that Theorem \ref{teo1} and Proposition \ref{teo2} are not specific of $\zeta$ per se: it's not hard for example to build from $\zeta$ other analytic functions which safisfy the functional equation and the property of being real on the real axis (and which might also have zeros not on the critical line), by multiplying it with an expression of the type $f(z)f(1-z)\overline{f(\bar{z})}\overline{f(1-\bar{z})}$. As \cite{sa03} and \cite{sa04} say, the Riemann Hypothesis would follow if we could improve these results and obtain a strict inequality
($|\zeta(\sigma+it)|>|\zeta(1-\sigma+it)|$ for $0<\sigma<1/2$ and $t>2\pi+1$), since roots have to be symmetric with respect to the critical line.

\end{rmk}

Let us see now how Theorem \ref{teo1} can bring to other formulations of the Riemann Hypothesis.

We suppose that the multiplicity of $s_0$ is $m=m(s_0)$ (then $m$ is as well the multiplicity of $1-s_0$, as we already remarked in the previous section). Since
$\zeta(s)$ and $\zeta(1-s)$ are analytic at $s_0$, the following
power series expansions hold for $s$ in a neighboorhood of $s_0$:
$$
\zeta(s)=\sum_{n=0}^{\infty}\frac{\zeta^{(n)}(s_0)}{n!}(s-s_0)^n,
$$
$$
\zeta(1-s)=\sum_{n=0}^{\infty}\frac{(-1)^n\zeta^{(n)}(1-s_0)}{n!}(s-s_0)^n.
$$
Therefore
$$
|F(s_0)|=\left|\lim_{s \to s_0} \frac{\zeta(s)}{\zeta(1-s)}\right|= \left|\frac{\zeta^{(m)}(s_0)}{\zeta^{(m)}(1-s_0)}\right|.
$$
By $(6)$, $\overline{\zeta(s)}=\zeta(\bar{s})$ for any $s$, so
that $|\zeta^{(m)}(s)|=|\zeta^{(m)}(\bar{s})|$.


We then also have
$$
|F(s_0)|= \left|\frac{\zeta^{(m)}(s_0)}{\zeta^{(m)}(1-\bar{s_0})}\right|.
$$


Thus we can write the following propositions:
\begin{prop}
The Riemann Hypothesis is true if and only if $\left|\frac{\zeta^{(m)}(s_0)}{\zeta^{(m)}(1-\bar{s_0})}\right|=1$ for any zero $s_0$ of $\zeta$.
\end{prop}
and 
\begin{prop}
The Riemann Hypothesis is true if and only if $\left|\frac{\zeta^{(m)}(s_0)}{\zeta^{(m)}(1-s_0)}\right|=1$ for any zero $s_0$ of $\zeta$.
\end{prop}

We can write explicit expansions for these derivatives. For example we can look at $F(s)$ with the representation of $\zeta$ as in
$(5)$. Then
$F(s)=\frac{\phi(s)}{\phi(1-s)}\cdot\frac{(1-2^{s})}{(1-2^{1-s})}$,
where $F(s_0)$, as above, is clearly meant to be its analytic continuation
in $s_0$. 

Since there's only one possible analytic continuation
and since the second fraction is a nonzero constant at $s_0$, then
$\frac{\phi(s_0)}{\phi(1-s_0)}\doteq \lim_{s \to s_0} \frac{\phi(s)}{\phi(1-s)}$ should also be some nonzero
constant. 

And, similarly as above,
$$
\\\left|\lim_{s \to s_0} \frac{\phi(s)}{\phi(1-s)}\right|= \left|\frac{\phi^{(m)}(s_0)}{\phi^{(m)}(1-s_0)}\right|.
$$

By the analytic properties of Dirichlet series (see \cite[section 11.7]{ap76}), we can write (differentiating term by term) for $\Re(s)>0$:
$$
\phi^{(m)}(s)=\sum_{n=1}^{\infty}\frac{(-1)^{n+m-1}(\log n)^m}{n^s}.
$$

So we would like to evaluate:
$$
\lim_{N\to\infty} \left| \frac{\sum_{n=1}^{N}(-1)^{n+m-1}(\log
    n)^m n^{-s_0}}{\sum_{n=1}^{N}(-1)^{n+m-1}(\log n)^m n^{-1+s_0}} \right|
$$

(limit of the denominator is nonzero, so that we can put the limit
before the fraction).

We can then rephrase the previous proposition in the following way:
\begin{prop}
Suppose that
$$
\sum_{n=1}^{\infty}(-1)^{n-1} n^{-s_0}=\sum_{n=1}^{\infty}(-1)^{n-1} n^{-1+s_0}=0
$$ 
with $s_0$ a zero of order $m$. Then this implies
$$
\lim_{N\to\infty} \left| \frac{\sum_{n=1}^{N}(-1)^{n+m-1}(\log
    n)^m n^{-s_0}}{\sum_{n=1}^{N}(-1)^{n+m-1}(\log n)^m n^{-1+s_0}} \right|=\left|\frac{(1-2^{1-s_0})}{(1-2^{s_0})}\right|
$$
if and only if the Riemann Hypothesis is true.
\end{prop}

Clearly a similar argument would hold if we consider $1-\bar{s_0}$ instead.

On the other hand it would be sufficient, in order to prove the Riemann Hypothesis, to show that $|F(s_0)|$ or
$$
\lim_{N\to\infty} \left| \frac{\sum_{n=1}^{N}(-1)^{n+m-1}(\log
    n)^m n^{-s_0}}{\sum_{n=1}^{N}(-1)^{n+m-1}(\log n)^m n^{-1+s_0}} \right|
$$
is going to be $0$ or $\infty$ (thus not a nonzero constant) if
$s_0$ is not on the critical line.




To investigate whether $|F(s_0)|$ is equal to $1$, we can consider also the function we built in the previous section through the approximants $H_N$.


For example showing that $H^{|F(s_0)|,|F(1-s_0)|}(s)$, which is
continuous in $s_0$, is also continuous at least in $1-s_0$,
would prove that $|F(s_0)|=1$, since we would have
$|F(1-s_0)|=|F(s_0)|$, because of continuity and Lemma \ref{lemma1}, but on the other hand $|F(s_0)|=1/|F(1-s_0)|$.

Or showing that there exists
a constant $k$ such that $H^{|F(s)|,|F(s)|}(s)=H^{k,k}(s)$
for all $s$ in the strip would also imply $|F(s_0)|=1$. (Clearly this holds locally for a neighboorhood: being $s$ in a sufficiently
small neighboorhood of $s_*$ ensures us that, if $s_*$ is a zero,
then by $(4)$, $\zeta(s)$ will not be zero, and, if $s_*$ is not
a zero, then by continuity $\zeta(s)$ will not be zero either. So
we are in the same situation as in the previous propositions
\ref{prop1} and \ref{prop2}, where what essentially counts is that the $H$-subscript
$N$ or $f(N,s)$ tends to infinity with $N$.)

In other words in this case Theorem \ref{teo1} can be reformulated as
\begin{prop}
The Riemann Hypothesis is true if and only if $H^{1,1}(s)$ is continuous throughout the strip.
\end{prop}
Notice that we can also state the following, which though is essentially weaker, as the continuity in this case appears to be much less at hand.
\begin{prop}
The Riemann Hypothesis is true if and only if $\lim_{N\to\infty} \left|\frac{H_N(s)}{H_N(1-s)}\right|$ is continuous throughout the strip.
\end{prop}
\begin{proof}
If the Riemann Hypothesis is true, then there are only zeros on the critical line, where $\lim_{N\to\infty} \left|\frac{H_N(s)}{H_N(\bar{s})}\right|=1$ since, for each $N$, $H_N(s)$ is the conjugate of $H_N(\bar{s})$ and so it has the same modulus. On the other side, we show that, if the Riemann Hypothesis were not true, then the limit would not be continuous. This is because in a zero out of the critical line, if we compute similarly as done in Lemma \ref{lemma1}, the limit would now be $0$ or $\infty$, while we know that $F(s)$ is always nonzero.
\end{proof}


\section{Acknowledgements}

I would like to thank everybody who has somehow helped me, even if I won't mention him personally, from the colleagues with whom I had fruitful discussions to the people who read part of or the whole manuscript, those who introduced me to the subject and everyone who has encouraged or supported me. Many thanks, among the others, to: Joachim Rosenthal, Alessandro Cobbe, Jens Zumbr\"agel, Felix Fontein, G\'erard Maze, Alessio Martini, Emanuele Spadaro, Ivan Contreras, Elisa Gorla, Patrick Sol\'e, Ashkan Nikeghbali, Christopher Hughes, Michele Elia, Pietro Peterlongo, Kim Dang.





\bibliography{huge} \bibliographystyle{plain}
\end{document}